\documentclass{ifacconf}

\usepackage{graphicx}      
\usepackage{natbib}        
\usepackage{amsmath}
\usepackage{amsfonts}
\usepackage{xcolor}

\usepackage{booktabs}
\usepackage{multirow}

\usepackage{siunitx}  

\usepackage{tikz}
\usetikzlibrary{external}
\usepackage{pgfplots}
\pgfplotsset{compat=1.16}
\usetikzlibrary{patterns}
\usepgfplotslibrary{fillbetween}

\usepackage{currfile}

\newcommand{\nthpoint}{20}
\begin{document}
\begin{frontmatter}

\title{Observer Design for an Inertia Wheel Pendulum with Static Friction\thanksref{footnoteinfo}} 

\thanks[footnoteinfo]{This work has been supported by the COMET-K2 Center of the Linz Center of Mechatronics (LCM) funded by the Austrian federal government and the federal state of Upper Austria.
}

\author[First]{L. Ecker} 
\author[First]{K. Schlacher}
\author[First]{M. Schöberl} 

\address[First]{Institute of Automatic Control and Control Systems Technology, Johannes Kepler University Linz,  Altenberger Str. 69, A-4040 Linz, Austria (e-mail: \{lukas.ecker, kurt.schlacher, markus.schoeberl\}@jku.at).}

\begin{abstract} 
A state observer design for the inertia wheel pendulum considering static friction of the actuated inertia disc is presented.
The frictional force is modeled by the Stribeck effect, with two separate differential equations describing the sticking and non-sticking system, respectively.
The transition between the two scenarios is in general determined by a static friction condition.
Three proposed observers are designed for both the sticking and non-sticking model.
Depending on the probability of the respective model given a current measurement, the more likely is selected for state estimation.
The performance is demonstrated on a laboratory demonstrator, where the observers including the proposed probabilistic model selection are compared.
\end{abstract}

\begin{keyword}
 Nonlinear Observer Design, Static Friction, State Observation, Inertia Wheel Pendulum, Bayesian Statistics 
\end{keyword}

\end{frontmatter}
\section{Introduction}

State estimation is a common method for determining non-measurable and unknown system states.
Just as the state estimation problem is generally more difficult in the nonlinear scenario than in the linear scenario, the same is true for the distinction between infinite- and finite-dimensional systems, see \cite{Meurer:2005}.
Moreover, in particular systems governed by a family of differential equations with algebraic switching conditions, i.\,e., switched or hybrid systems, see, e.\,g., \cite{Liberzon:2003}, require special attention. 
Switching conditions, as encountered when modeling mechanical systems with static friction, are generally state-dependent.
Thus, depending on the current values of the state, the system behaves according to a different differential equation. 
This is usually not a problem when evaluating the initial value problem, since the state variables required for the switching condition are known from the simulation. 
However, for the state estimation problem, the state-dependent condition is critical because, except for certain measured quantities, only the estimated quantities are known to the observer. Depending on which differential equation is used, different estimates result, which in turn have different implications for the choice of the actual model for the observer.
The present paper focuses on the observer design for the inertia wheel pendulum, a nonlinear mechanical benchmark example in control theory as found in, e.\,g., \cite{Spong:1999}, with switching sticking and non-sticking subsystem.
The stiction of the wheel, which typically occurs during pendulum tests in real experimental setups, is modeled as an alternation between two models, with each described by a separate differential equation. 
An extended Kalman filter (EKF), a nonlinear observer with linear error dynamics ($\text{NO}_\text{L}$), and a nonlinear observer with passive error dynamics ($\text{NO}_\text{P}$), are designed and evaluated on the laboratory demonstrator.
The widely used EKF, see, e.\,g., \cite{Simon:06}, is based on a sampled-data system representation and stochastic methods, while the $\text{NO}_\text{L}$ and $\text{NO}_\text{P}$ exploit system theoretic approaches such as linearization by output injection, see, e.\,g., \cite{KRENER:1983}, and passivity properties,  see, e.\,g., \cite{Schaft:2017}. 
Therefore, the observer with passive error dynamics exploits a port-Hamiltonian representation of the inertia wheel pendulum, see, e.\,g.,\cite{Ortega:2001}, which allows to inject a dissipative behavior of the error system between state and estimate. 
The port-Hamiltonian framework has proven useful for both finite and infinite dimensions in system analysis and in the design of controllers and observers, see, e.\,g., \cite{Macchelli:2004}, \cite{SCHOBERL:2014}, \cite{Malzer:2020}, and \cite{Ecker:2021}.

The paper considers a finite-dimensional switched system.
The main difficulties that arise in the transition between the sticking and non-sticking model for the observer are addressed by probabilistic model selection. Instead of the sole evaluation of the switching condition with estimates from either the sticking or the non-sticking model, the marginal likelihood of the two models given the obtained measurement is used. The Bayesian factor, a ratio of the likelihoods of the competing models, see, e.\,g., {\cite{Bernardo:2002}}, determines which model is used for the prediction of the state.
The paper is structured as follows. In Section 2, the model of the inertia wheel pendulum including static friction of the wheel is briefly derived. Section 3, is dedicated to the design of the observers as well as the presentation of the probabilistic model selection. Finally, in Section 4, the performance of the proposed observers are demonstrated on a laboratory model.


\section{Mathematical modeling}

This section is dedicated to the mathematical modeling of the inertia wheel pendulum. Since the mechanical system is a common example in control theory, see, e.\,g., 
\cite{Spong:1999}, a detailed derivation of the equation of motion is omitted and more attention is devoted to the modeling of the static friction.

\subsection{Inertia Wheel Pendulum}
\begin{figure}
	\begin{center}
		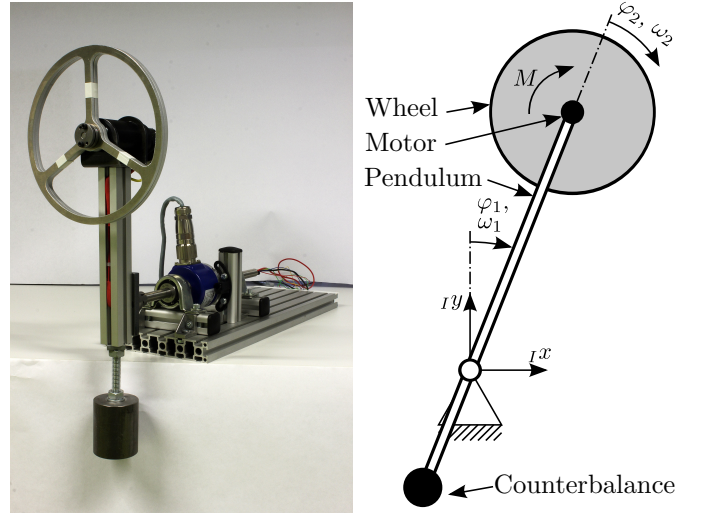
		\caption{{Schematic of the inertia wheel pendulum.}}
		\label{fig:schematic_and_image_inertia_wheel_pendulum_en}
	\end{center}
\end{figure}

The model of the inertia wheel pendulum under consideration, see, Fig. \ref{fig:schematic_and_image_inertia_wheel_pendulum_en}, consists of a rigid, ball-bearing pendulum with actuated wheel at one end and counterweight at the opposite end. 
The absolute angle of the pendulum is denoted by $\varphi_1$ and the relative angle of the wheel by $\varphi_2$.
The wheel is driven by a subordinate torque-controlled DC electric motor.
The electrical dynamics, which are much faster compared to the dynamics of the mechanical system, are neglected.
Therefore, the motor torque $M$ is considered as the input of the model of the inertia wheel pendulum.
The governing equations of motion with regard to the generalized coordinates $q = [q^\alpha] = [\varphi_1,\varphi_2]^T$ and velocities $\dot{q} = [\dot{q}^\alpha] =  [\dot{\varphi_1},\dot{\varphi_2}]^T$ are derived by means of the Euler-Lagrange equations
\begin{equation}
\label{equ:lagrange_equation}
\frac{\mathrm{d}}{\mathrm{d}t}\frac{\partial L}{\partial \dot{q}^\alpha} - \frac{\partial L}{\partial q^\alpha} = Q_\alpha, \quad \alpha = 1,2\,,
\end{equation}
where the Lagrange function $L = T - V$ is the difference between the kinetic energy $T = \frac{1}{2} \theta_1 \dot{\varphi}_1^2 + \frac{1}{2} \theta_2 \left(\dot{\varphi}_1 + \dot{\varphi}_2\right)^2$ and the potential energy $V = a \cos(\varphi_1)$.
The moments of inertia of the subsystem pendulum are summarized by the parameter $\theta_{1}$, whereas the moment of inertia of the wheel with respect to the center of gravity is described by $\theta_2$.
Moreover, the factor $a$ covers the components of the potential energy  of the overall system. 
The viscous damping between the pendulum and the base as well as the pendulum and the inertia wheel are taken into account.  
The damping terms are assumed to be proportional to the respective angular velocity, i.\,e., $d_\alpha \dot{\varphi}_\alpha $, which together with the input moment $M$ and friction moment $M_S$ form the terms of the generalized forces $Q = [Q_\alpha] = [ -d_1\dot{\varphi}_1, -d_2 \dot{\varphi}_2 + M {- M_S}]^T$.
Evaluating the Euler-Lagrange equations (\ref{equ:lagrange_equation}), the equations of motion correspond to the second order differential equation system
\begin{equation}\small
\begin{bmatrix}
\theta_1 + \theta_2 & \theta_2\\
\theta_2 & \theta_2
\end{bmatrix}
\begin{bmatrix}
\ddot{\varphi}_{1}\\
\ddot{\varphi}_{2}
\end{bmatrix}
-
\begin{bmatrix}
a \sin(\varphi_1)\\
0
\end{bmatrix}
=
\begin{bmatrix}
-d_1 \dot{\varphi_1} \\ -d_2 \dot{\varphi_2} + M {- M_S}
\end{bmatrix}.
\label{equ:iwp_eom}
\end{equation}
Furthermore, the equations can be rewritten as a system of first order differential equations with state $x = [\varphi_{1} ,\varphi_{2} , \omega_1 , \omega_2 ]^T$, where $\omega_\alpha = \dot{\varphi}_\alpha$. Since the resulting system does not depend on the angle $\varphi_2$, the  equation $\dot{\varphi}_2 = \omega_2$ is omitted and only the three dimensional system 
\begin{equation}\hspace*{-0.2cm}
\label{equ:non_adherent_system}
\begin{bmatrix}
\dot{\varphi}_{1}\\
\dot{\omega}_{1}\\
\dot{\omega}_{2}
\end{bmatrix}
=
\begin{bmatrix}
\omega_{1}\\
\frac{1}{\theta_{1}} \left( a\sin \left( \varphi_{1} \right) -d_{1}\omega_{1} +d_{2}\omega_{2}-M {+ M_S}\right)\\
-\frac{a}{\theta_{1}}\sin(\varphi_{1}) + \frac{d_1}{\theta_{1}}  \omega_{1} +\frac{1}{\theta_{c}} (M {- M_S} - d_2 \omega_{2})
\end{bmatrix}
\end{equation}
with state $x = [\varphi_{1}, \omega_1, \omega_2 ]^T$ and abbreviation $\theta_c = \frac{\theta_{1} \theta_{2}}{\theta_{1}+\theta_{2}}$ is considered.

\subsection{Static friction}

The present model is not yet satisfactory for the state observation of the inertia wheel pendulum.
The occurring discrepancy between the mathematical model and real measured data, which become apparent during the validation (see Section \ref{sec:validation}), indicate that crucial effects are not taken into account by the model (\ref{equ:non_adherent_system}).
Depending on the trajectory of the system, i.\,e., usually the drop down or swing up of the pendulum, different disturbances are crucial. In the swing-up experiment, the static friction of the pendulum {dominates}, while in the drop-down experiment, the static friction of the wheel {dominates}. However, the present paper focuses on the observer design with the pendulum angle as measured quantity. Since a direct measurement of the wheel angle is not considered, the static friction of the wheel is more critical for the observer design and therefore the main emphasis of this work. 
The modeling of the static friction of the inertia wheel is realized by a structural {switch between two differential equations.}
The friction force is described as a function of the friction velocity by the Stribeck curve
\begin{align*}
M_S(\omega_2) = r_C \,\mathrm{sgn}(\omega_2) + (r_S-r_C) e^{-\left(\frac{\omega_2}{\omega_{2,0}}\right)^2}\mathrm{sgn}(\omega_2)
\end{align*}
where $r_C$ corresponds to the Coulomb friction and $r_S$ to the Stiction parameter. The negative slope friction for a small transition of velocity from zero is referred to as the Stribeck effect, see \cite{Ellis:2003}. Note that viscous damping has already been taken into account in the derivation of the mathematical model. Static friction, i.\,e., $\omega_2 = 0$ and $\dot{\omega}_2 = 0$, must meet the condition
\begin{align}\label{equ:stiction_condition}{\left|\frac{\theta_2}{\theta_1 +\theta2}\left(-a \sin(\varphi_{1})+d_1\omega_1\right)+M\right| < r_S.}
\end{align}
The model of the inertia wheel pendulum in stiction with $\omega_{2} = 0$ and {$M_S = 0$} is given by 
\begin{equation}\label{equ:adherent_system}
\begin{bmatrix}
\dot{\varphi}_{1}\\
\dot{\omega}_{1}\\
\dot{\omega}_{2}
\end{bmatrix}
=
\begin{bmatrix}
\omega_{1}\\
{\frac{1}{\theta_{1}} \left( a\sin \left( \varphi_{1} \right) -d_{1}\omega_{1} \right)}\\
0
\end{bmatrix},
\end{equation}
where the associated Lagrange function is equivalent to
\begin{equation}
\label{equ:stiction_lagrangian}
L = \frac{1}{2}\theta_{1}\dot{\varphi}_1^2 - a\cos(\varphi_{1}).
\end{equation}

\subsection{Port-Hamiltonian Representation}

The subsequent nonlinear observer design with passive error dynamics will exploit a port-Hamiltonian representation of the mechanical system. By means of the regular Legendre mapping the Hamiltonian $H = \sum_{\alpha = 1}^{n}p_\alpha \dot{q}^\alpha - L$ of the non-sticking system (\ref{equ:non_adherent_system}) corresponds to
\begin{equation}
H = \frac{1}{2}\frac{1}{\theta_{1}}\left(p_1-p_2\right)^2 + \frac{1}{2}\frac{1}{\theta_{2}}p_2^2 + a \cos(\varphi_{1})
\end{equation}
with the generalized momenta $p_\alpha = \frac{\partial L}{\partial \dot{q}^\alpha}$ introduced as
\begin{align}
	\begin{bmatrix}
p_{1}\\
p_{2}
\end{bmatrix}
= 
\begin{bmatrix}
\theta_1 + \theta_2 & \theta_2 \\
\theta_2 & \theta_2
\end{bmatrix}
\begin{bmatrix}
\dot{\varphi}_{1}\\
\dot{\varphi}_{2}
\end{bmatrix}.
\end{align} 
The Hamiltonian formalism $\dot{p}_\alpha = \frac{\partial L }{\partial q^\alpha} + Q_\alpha$ allows to represent the pendulum as port-Hamiltonian system 
\begin{equation}  \label{equ:ph_non_sticking}
\dot{x} = \left(J - R\right) \left(\frac{\partial H}{\partial x}\right)^T + G u
\end{equation}
with state {$x = [{\varphi_{1}},{p}_1, {p}_2]^T$}, {input $u = M - M_S$}, partial derivatives {${\partial_x H} = {\partial/ \partial_x H}= [-a \sin(\varphi_1),
(p_1-p_2)/{\theta_1},
-(p_1-p_2)/{\theta_1} + p_2/\theta_{2}
]$} and matrices
\begin{equation*}
	{J=-J^T} = 
	\begin{bmatrix}
	0 & 1 & 0 \\
	-1 & 0 & 0 \\
	0 & 0 & 0
	\end{bmatrix},\;
{R=R^T} =
\begin{bmatrix}
	0 & 0 & 0 \\
	0 & d_1 & 0 \\
	0 & 0 & d_2
	\end{bmatrix}
,\;
G = 
\begin{bmatrix} 0 \\ 0\\ 1	\end{bmatrix}.
\end{equation*}
Proceeding similarly, the port-Hamiltonian representation of the inertia wheel pendulum in stiction with Lagrangian (\ref{equ:stiction_lagrangian}), generalized momentum $p_1 = \frac{\partial L}{\partial \dot{\varphi}_1} =  \theta_{1} \dot{\varphi}^1$ and Hamiltonian $H = \frac{1}{2}\frac{1}{\theta_{1}}\dot{p}_{1}^2 + a\cos(\varphi_{1})$ is given by
\begin{equation} \label{equ:ph_sticking}
\underbrace{\begin{bmatrix}
	\dot{\varphi}_1\\
	\dot{p}_{1}
	\end{bmatrix}}_{\dot{x}}
=
{\biggl(
\underbrace{	\begin{bmatrix}
	0 & 1 \\
	-1 & 0 
	\end{bmatrix}}_{J=-J^T}}
-
\underbrace{\begin{bmatrix}
	d_1 & 0 \\
	0 & 0
	\end{bmatrix}}_{R=R^T}
\biggl)\cdot
\underbrace{	\biggl(
	\begin{matrix}
	-a \sin(\varphi_1)\\
	\frac{1}{\theta_1}p_1
	
	\end{matrix}
	\biggl)}_{\left(\frac{\partial H}{\partial x}\right)^T}.
\end{equation}


\section{Observer design}\label{sec:observer}
 
The presented observer design focuses on the estimation of the state $x = [\varphi_{1},\omega_{1},\omega_{2}]^T$ with the pendulum angle $\varphi_{1}$ as measured output. For this purpose, three nonlinear state estimators for the inertia wheel pendulum are introduced. 

\subsection{Extended Kalman Filter}\label{sec:EKF}

The Kalman filter is an optimal state observer for linear systems with Gaussian noise. Under the assumption that both the probability density function of the initial state and  the noise are Gaussian distributed, the first two stochastic moments, i.\,e., the expectation and the variance, are sufficient for a stochastic representation. 
The extended Kalman filter is a nonlinear version of the Kalman filter based on the linearization of the dynamics about the current estimate. Since the EKF is de facto standard in the theory of nonlinear state estimation, a detailed derivation, as found in \cite{Simon:06}, is omitted.
The time-continuous models (\ref{equ:non_adherent_system}) and (\ref{equ:adherent_system}) are discretized using a fourth-order Runge-Kutta schema with sampling time $T_a = 5\,\mathrm{ms}$, i.\,e., considered are the equations
\begin{subequations}
	\label{equ:nonlinear_discrete_systems}
	\begin{align}
	x_{k+1} &= f_k(x_k,u_k,w_k) \\
	y_k &= h_k(x_k,v_k)
	\end{align}
\end{subequations}
with time index $k$, system state  
$x_k = [\varphi_{1},\omega_{1},\omega_{2}]^T$, output $y_k  = \varphi_{1}$, input $u_k = M $ as well as process noise $w_k $ and measurement noise $v_k$. 
The estimation of the state is denoted with $\hat{x}_k$. 
Depending on the number of measurements included, it is referred to as an a posteriori or a priori estimate. The a priori state $\hat{x}_k^-$ considers all measurements up to the previous time $k-1$, i.\,e., $Y_{k-1} = \{y_0,\ldots,y_{k-1}\}$, where as the a posteriori state $\hat{x}_k^+$ includes all measurements up to actual time $k$, i.\,e., $Y_{k} = \{y_0,\ldots,y_{k}\}$. The Taylor series expansions of the sampled-data system (\ref{equ:nonlinear_discrete_systems}) along the estimates corresponds to the linearised system
	\begin{subequations}
		\begin{align}
		x_{k+1} &= F_k x_k + L_k w_k + \tilde{u}_k \\
		y_{k} &= H_k x_k + M_k v_k 
		\end{align}
	\end{subequations}
	with Jacobian matrices $ F_k = \left. \frac{\partial f_k}{\partial x}  \right \vert_{\hat{x}^+_k}$,  $L_k = \left. \frac{\partial f_k}{\partial w}  \right \vert_{\hat{x}^+_k}$ $
	H_k = \left. \frac{\partial h_k}{\partial x}  \right \vert_{\hat{x}^-_k}$, $ 
	M_k = \left. \frac{\partial h_k}{\partial v}  \right \vert_{\hat{x}^-_k}$
	and input $\tilde{u}_k = f_k(\hat{x}_k^+, u_k, 0) - F_k \hat{x}^+_k$.
Assuming mean-free Gaussian process and measurement noise, i.\,e., $ w_k \sim (0,Q_k)$ and $v_k \sim (0,R_k)$, the prediction step is defined as
$P^-_{k+1} = F_k P_k^+ F_k^T + L_k Q_k L_k^T$ and $\hat{x}^-_{k+1} = f_k(\hat{x}_k,u_k,0)$
while the measurement update can be expressed by
\begin{subequations}
	\begin{align}
	K_k &= P^-_k H_k^T(H_k P^-_k H_k^T + M_k R_k M_k^T)^{-1}\\
	\hat{x}^+_{k} &=  \hat{x}^-_k + K_k \left(y_k - h_k(\hat{x}^-_k,0)\right)\\
	P^+_k &= P^-_k - K_k H_k P^-_k.
	\end{align}
\end{subequations}


\subsection{Nonlinear Observer with linear error dynamics}

The nonlinear observer with linear error dynamics is based on the idea of a coordinate change that splits the system into a linear and a nonlinear term, where the nonlinear part may only depend on input and output. This linearization by output injection, see \cite{KRENER:1983} and \cite{Krener:1985}, has the advantage that if the observer is chosen as a copy of the plant in adapted coordinates, the nonlinear terms cancel. The nonlinear time-continuous models (\ref{equ:non_adherent_system}) and (\ref{equ:adherent_system}) can be represented by a regular coordinate transformation $z = \Phi(x)$ and output $y$ as 
$
\dot{z} = A z + \alpha(u,y),
$ 
where $\alpha(u,y)$ corresponds to the nonlinear term that depends only on the input $u$ and output $y$. In the case of the inertia wheel pendulum, the nonlinearities occur exclusively in the measurable state $y = \varphi_{1}$, which is why a desired splitting is already possible with $z = \Phi(x) = x$. Therefore, the decomposition of the inertia wheel pendulum in the non-sticking scenario with original coordinates is given by
{
\begin{equation*}
\underbrace{
	\begin{bmatrix}
	\dot{\varphi}_1\\
	\dot{\omega}_1 \\
	\dot{\omega}_2 
	\end{bmatrix}}_{\dot{x}}
= 
\underbrace{\renewcommand{\arraystretch}{1.2}
	\begin{bmatrix}
	0& 1 & 0\\
	0&-\frac {d_{1} }{\theta_1}& \frac {d_{2} }{\theta_1} \\
	0&\frac{d_1}{\theta_1}&-\frac{d_2}{\theta_c} \\
	\end{bmatrix}}_{A}
\underbrace{
	\begin{bmatrix}
	\varphi_1\\
	\omega_1 \\
	\omega_2 \\
	\end{bmatrix}}_{x}
+
\underbrace{\renewcommand{\arraystretch}{1.2}
	\begin{bmatrix}
	0 \\
	\frac{a}{\theta_{1}}\sin(\varphi_1) - \frac{M-M_S}{\theta_{1}} \\
	-\frac{a}{\theta_{1}}\sin(\varphi_1) + \frac{M-M_S}{\theta_{c}} 
	\end{bmatrix}}_{\alpha(u,y)},
\end{equation*}}
where the decomposition of the sticking model is given by
{\begin{equation*}
\underbrace{
	\begin{bmatrix}
	\dot{\varphi}_1\\
	\dot{\omega}_1 
	\end{bmatrix}}_{\dot{x}}
= 
\underbrace{\renewcommand{\arraystretch}{1.2}
	\begin{bmatrix}
	0& 1 \\
	0&-\frac {d_{1} }{\theta_1}  \\
	\end{bmatrix}}_{A}
\underbrace{
	\begin{bmatrix}
	\varphi_1\\
	\omega_1 \\
	\end{bmatrix}}_{x}
+
\underbrace{\renewcommand{\arraystretch}{1.2}
	\begin{bmatrix}
	0 \\
	\frac{a}{\theta_{1}}\sin(\varphi_1)  \\
	\end{bmatrix}}_{\alpha(u,y)}.
\end{equation*}}
The nonlinear observer with state $\hat{x}$ is composed of a copy of the plant and linear correction term
\begin{equation}\label{equ:nlowled}
\dot{\hat{x}} = A\hat{x} + \alpha(u,y) + K(y-\hat{y})
\end{equation}
where $\hat{y}$ corresponds to the estimated linear output $\hat{y} = C \hat{x}$  and $K$ represents the correction coefficients of the observer.
The dynamics of the linear error $e = x - \hat{x}$ are thus given by
\begin{equation}\label{equ:continious_linear_error_dynamics}
\dot{e} = (A-KC)e.
\end{equation}
Both the linear subsystem of the sticking and non-sticking models are observable, allowing observer design methods from linear control theory such as pole placement or optimal designs based on the solution of the algebraic Riccati equation, see, e.\,g., \cite{Simon:06}.
{The proposed nonlinear observer with linear error dynamics (\ref{equ:nlowled}) is implemented quasi-continuous.
The previously described design approach can also be carried out for the sampled-data system $x_{k+1} = A_dx_k + G_d \alpha(u_k,y_k)$  with $A_d = e^{A T_a}$ and $G_d = \int_{0}^{T_a}e^{A \tau}\mathrm{d}{\tau}$, where $\alpha(u_k,y_k)$ is assumed to be constant in the time interval $t_k$ and $t_{k+1}$. 
Note that since the state and hence the output are quantities of a continuous process, the term $\alpha(u,y)$ is generally not constant over the sampling period.
However, given sufficient sampling time, satisfactory observer performance can be achieved even with this non-exact approximation to the sampled-data system.
Nevertheless, the nonlinear discrete observer given by
\begin{equation}\label{equ:nlowled_sampled_data}
\hat{x}_{k+1} = A_d \hat{x}_k + G_d \alpha(u_k,y_k)+ K_d(y_k-\hat{y}_k)
\end{equation}
with correction term $K_d$ also leads to an error dynamics
\begin{equation}\label{equ:discrete_linear_error_dynamics}
e_{k+1} = (A_d-K_d C)e_k
\end{equation}
with linear behavior. Both the linear subsystem of the sticking and of the non-sticking model are observable.
Finally, the nonlinear observer is designed by an appropriate choice of $K$ or $K_d$, with the objective of stabilizing the linear error dynamics (\ref{equ:continious_linear_error_dynamics}) and (\ref{equ:discrete_linear_error_dynamics}), respectively.}

%

\subsection{Nonlinear Observer with passive error dynamics}
 
The underlying idea is to design a nonlinear observer that induces a passive behavior of the error of state and estimation. Exploiting the port-Hamiltonian representation of the system, the passive error dynamics is derived with the observer as a copy of the plant. In addition, a compensation term allows the imposition of an appropriate desired Hamiltonian function {$H_d$} {of the error system}, and a measurable collocated output can be used for the injection of damping behavior. In this way, the asymptotic stability of the error system {for the model in stiction and the model not in stiction} should be ensured. Considering the port-Hamiltonian representation  (\ref{equ:ph_non_sticking}) and (\ref{equ:ph_sticking}), the nonlinear observer with state $\hat{x}$ and system output $y = h(x)$ is introduced as
\begin{align}\label{equ:ph_observer}
	\dot{\hat{x}} &= \left(J-R\right)\left(\partial_{\hat{x}} H(\hat{x}) + \Phi(\hat{x},y)\right)^T + G u - G_o u_o,
\end{align}
where $\Phi(\hat{x},y)$ corresponds to the compensation term
and $G_o u_o$  to the observer input used for the damping injection. {Since the compensation term $\Phi(\hat{x},y)$ and the input $G_o u_o$ are quantities of the observer, they should only depend on the estimation $\hat{x}$ and the measurement $y$. As mentioned, the compensation term will be used to induce the desired Hamiltonian function $H_d$. For this purpose the compensation term is determined as \begin{equation*}
	\Phi = \left(\partial_x H(x)-\partial_{\hat{x}} H(\hat{x})- \partial_e H_d(e)\right)^T,
	\end{equation*} where $ e = x - \hat{x}$ corresponds to the considered linear error between state and estimate.
	As a result of the proposed observer (\ref{equ:ph_observer}) and the compensation term, the error dynamic is given by 
	\begin{align*}
	\dot{{e}} &= \left(J-R\right)\underbrace{\left(\partial_x H(x)-\partial_{\hat{x}} H(\hat{x})- \Phi(\hat{x},y)\right)^T}_{\left(\partial_e H_d(e)\right)^T}  + G_o u_o,
	\end{align*}
	where the evolution of the Hamiltonian corresponds to
	\begin{align*}
	\dot{H_d}(e) = -\left(\partial_e H_d(e)\right) R \left(\partial_e H_d(e)\right)^T + \left(\partial_e H_d(e)\right) G_o u_o.
	\end{align*}For the implementation of the observer, an appropriate positive definite Hamilton function  $H_d>0$ must ensure that the compensation term depends only on the estimate $\hat{x}$ and the output $y$, i.\,e., $\Phi(\hat{x},y)$ holds.
	This functional dependency can be verified by means of exterior algebra, where $\mathrm{d}$ refers to the exterior derivative with respect to the coordinates $x$ and $\wedge$ to the exterior product. 
	Therefore, considering in general $m$ output equations, the entries of the compensation term $\Phi^i$ defined by $H_d$ with $1\leq i \leq n$ must satisfy the condition
	$
	\mathrm{d}\Phi^i \wedge \mathrm{d} h^1 \wedge \ldots \wedge \mathrm{d}h^m = 0.
	$
	In order to verify the asymptotic stability of the error system, not only the positive definiteness of $H_d$ but also the negative definiteness of $\dot{H}_d$ must be ensured.
	Since for mechanical systems in general only $R \ge 0$ holds, see (\ref{equ:ph_non_sticking}) and (\ref{equ:ph_sticking}), and thus $\dot{H}_d \le 0$, additional dissipative behavior has to be included to force a steady decrease of the observation error.
	Therefore, a desired matrix $R_d = R_d^T = R + \bar{R} > 0$
	must match the condition
	\begin{multline*}
	\hspace*{-0.4cm} -\left(\partial_e H_d(e)\right) R \left(\partial_e H_d(e)\right)^T + \underbrace{\left(\partial_e H_d(e)\right) G_o u_o}_{F(\hat{x},u,y)} \\ {=} -\left(\partial_e H_d(e)\right) R_d \left(\partial_e H_d(e)\right)^T,
	\end{multline*}
	which allows the observer input with
	\begin{align*}
	G_o u_o = F(\hat{x},u,y) = -\bar{R}\left(\partial_e H_d(e)\right)^T
	\end{align*}
	asymptotic stabilize the error dynamics. Note that $G_0 u_0$ as part of the observer must also depend only on $\hat{x}$, $u$ and $y$.
	Considering the inertia wheel pendulum with measured output $y = h(x) = \varphi_{1}$, the positive definite function
	$
	H_d(e) = \frac{\alpha}{2}\tilde{\varphi}_{1}^2 + \frac{1}{2 \theta_1}\left(\tilde{p}_1-\tilde{p}_2\right)^2 + \frac{1}{2 \theta_2} \tilde{p}_2^2
	$
	with $\alpha > 0$ serves as proper desired Hamiltonian, since $\mathrm{d}\Phi^i \wedge \mathrm{d} h = 0$ holds. Likewise, the matching condition is fulfilled with
	\begin{align*}
	G_o u_o = -\bar{R}\left(\partial_e H_d(e)\right)^T = -\begin{bmatrix}
	\beta & 0 & 0\\
	0 & 0 & 0\\
	0 & 0 & 0
	\end{bmatrix}
	\left(
	\begin{matrix}
	\alpha \tilde{\varphi}_1\\
	. \\
	.
	\end{matrix}
	\right),
	\end{align*}
	where also simultaneously $R_d = R + \bar{R} > 0$ ensures the desired behavior of steady decreasing estimation error. The desired Hamiltonian for the model in stiction is set to
	$H_d(e) = \frac{\alpha}{2}\tilde{\varphi}_{1}^2 + \frac{1}{2 \theta_1} \tilde{p}_1^2$
	and the observer input to
	\begin{align*}
	G_o u_o = -\bar{R}\left(\partial_e H_d(e)\right)^T = -\begin{bmatrix}
	\beta &  0\\
	0 & 0 \\
	\end{bmatrix}
	\left(
	\begin{matrix}
	\alpha \tilde{\varphi}_1\\
	. \\
	\end{matrix}
	\right).
	\end{align*}
	The nonlinear observer with passive error dynamics (\ref{equ:ph_observer}) is implemented quasi-continuous.
}
\subsection{Model selection}
This section presents an approach for the implementation of a switching law between two models in the scope of a observer design with no complete knowledge of the state. Since the state dependence of the derived switching condition in the observer must be evaluated with estimates, a correct selection of the currently correct model for the observer is not always possible. Therefore, a probabilistic model selection based on the Bayesian factor is proposed. The basic idea is that with the current measurement obtained, the probability of both models is determined and the more likely one is used for prediction. The output is interpreted as a random variable $Y$ because a perturbation $v_k$ with probability density function $p_V(v_k)$ enters the output equation with $y_k = h(x_k) + v_k$. The non-sticking model (\ref{equ:non_adherent_system}) is referred to as $M_{1}$, and the sticking model (\ref{equ:adherent_system}) is referred to as $M_2$, where $f_{M_1,k}$ and $f_{M_2,k}$ denote their corresponding sampled-data systems. In principle, this methodology is not limited to two models. Nevertheless, considered are the two models $M_1$ and $M_2$ as well as the a posterior estimate of the previous time step $\hat{x}_{k-1}$ as well as the current measurement $y_k$. The marginal likelihoods of the two models given the obtained measurement in ratio results in the so-called Bayes factor
\begin{equation}
K = \frac{\mathrm{Pr}(M_1 | Y = y_k)}{\mathrm{Pr}(M_2 | Y = y_k)}.
\end{equation}
Therefore, when $K>1$, it is more likely that the measurement $y_k$ stems from the model $M_1$, and vice versa, when $K<1$, the model $M_2$ is to be preferred. With equal probable models, i.\,e., $\mathrm{Pr}(M_1) = \mathrm{Pr}(M_2)$ and Bayes' theorem, see, e.\,g., \cite{Papoulis:2002},
\begin{equation}
\mathrm{Pr}(M | Y = y_k) = \frac{p_Y(y_k|M)\mathrm{Pr}(M)}{p_Y(y_k)},
\end{equation}
where $p_Y(y_k) = p_Y(y_k|M_1)\mathrm{Pr}(M_1) + p_Y(y_k|M_2)\mathrm{Pr}(M_2)$ denotes the probability density function of the output, the ratio of the Bayes factor can be simplified to  
\begin{equation}
K = \frac{p_Y(y_k|M_1)}{p_Y(y_k|M_2)}.
\end{equation}
For the computational implementation the conditional probability density function can be expressed in terms of the probability density of the perturbation as 
\begin{equation}
p_Y(y_k|M_d) = p_V(v_k) \circ \left(y_k - h(\hat{x}_{M_d,k})\right)
\end{equation}
with $\hat{x}_{M_d,k} = f_{M_d,k}(\hat{x}_{k-1})$. 

\section{Experimental Results}\label{sec:experimental_results}

This section is dedicated to the validation of the proposed models and the presentation of the observer performance by the drop-down experiment of the inertia pendulum.

\subsection{Validation}\label{sec:validation}
The validation should highlight the discrepancy between the model without static friction and the real measured data in order to justify the additional effort for modeling static friction. Furthermore, the validation allowed to identify the parameters of the Stribeck curve as well as the viscous damping coefficients. The parameters of the static friction model as well as the experimental determined parameters of the inertia wheel pendulum are listed in Table \ref{tab:parameter}. Note that the stiction of the pendulum as well as the counter-torque due to the supply lines to the drive wheel are not taken into account. However, in the course of validation, the initialization error of the incremental encoder was  compensated by calibration.
\begin{table}[ht]
	\begin{center}
			
	\caption{Parameters of the inertia wheel pendulum and the static friction model.}\label{tab:parameter}
	{
	\begin{tabular}{llrcllr}
		\toprule
		
		 \multicolumn{3}{l}{Model parameters} && \multicolumn{3}{l}{Static friction model} \\ 
		\cmidrule{1-3}
		\cmidrule{5-7}
		  \multicolumn{2}{l}{Parameter}  & Value  & &  \multicolumn{2}{l}{Parameter}  & Value \\
		\midrule 
		
		 $a$ & [\SI{}{\kilogram \square \metre \per \square \second}] &0.15535 &&  $r_C$&[\SI{}{\kilogram \square \metre \per \square \second}]& 0.0024 \\
		 $\theta_{1}$ & [\SI{}{\kilogram \square \metre }] & 0.05045 && $r_S$ &[\SI{}{\kilogram \square \metre \per \square \second}]&0.0026\\
		 $\theta_{2}$ & [\SI{}{\kilogram \square \metre }] & 0.00113 && $w_{2,0}$&[\SI{}{\per \second}]&0.0501 \\
		 $d_1$ & [\SI{}{\kilogram \square \metre \per \second}] & 0.00885 && & \\
		 $d_2$ & [\SI{}{\kilogram \square \metre \per \second}] & 0.00015 && & \\
		\bottomrule
	\end{tabular}}
	\end{center}

\end{table}
The reference measurement data, i.\,e., the angles of pendulum and wheel, of the real laboratory model were recorded by two incremental encoders. The associated angular velocities were determined by gate time measurement.
The following Figure \ref{fig:validation} shows the validation of both the model without static friction and the model with static friction and structure switch based on (\ref{equ:stiction_condition}). {While the results of measurement and simulation of both models hardly differ for $\varphi_{1}$ and $\omega_{1}$, significant discrepancies between measurement and the model without static friction occur for $\omega_{2}$.
In contrast, the model with static friction reproduces the real behavior of the pendulum very well. Therefore, it is also used for the observer design, since the estimation performance depends to a high degree on the model accuracy.}

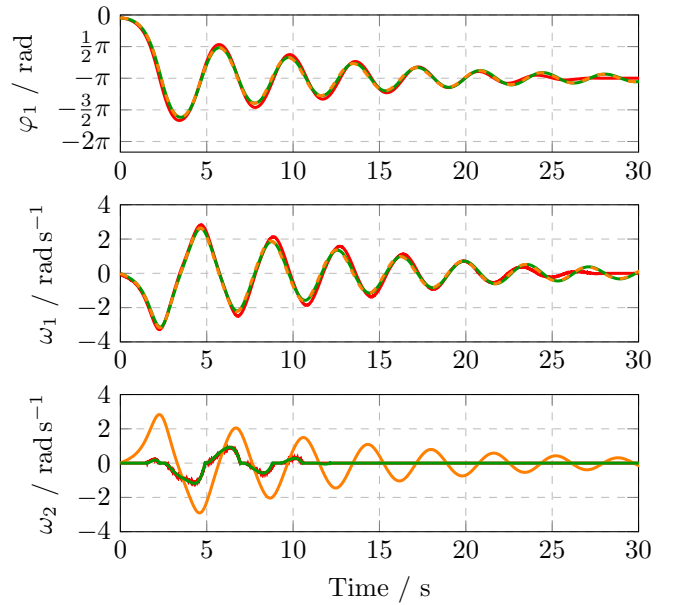
\begin{figure}[ht]
	\begin{center}
\tikzsetnextfilename{\currfilebase}
\begin{tikzpicture}

	\begin{axis}[
		name = one,
		width=8.4cm,
		height=3.4cm,
		ylabel={$\varphi_{1}$ / \SI{}{\radian}},
		xmin=0, xmax=30,
		ymin=-6.81, ymax=0,
		xtick={0,5,10,15,20,25,30},
		ytick = {-2*pi,-1.5*pi,-pi,-0.5*pi,0},
		yticklabels = {$- 2\pi$, $-\frac{3}{2}\pi$,$- \pi$,$\frac{1}{2}\pi$, $ 0$},
		legend pos=north west,
		ymajorgrids=true,
		xmajorgrids=true,
		grid style=dashed,
		anchor=above north
	]
	
	\addplot[color=red,	mark=none, line width=0.4mm] table[x=time, y=data (phi1), col sep=comma] {./data/data_meas.csv};
	
	\addplot[color=green!60!black, mark=none, line width=0.4mm] table[x=time, y=data (phi1), col sep=comma] {./data/data_model_stiction.csv};
	
	\addplot[color=orange, dashed, mark=none, line width=0.4mm] table[x=time, y=data (phi1), col sep=comma] {./data/data_model.csv};
	
	\end{axis}

	\begin{axis}[
		name = two,
		width=8.4cm,
		height=3.4cm,
		ylabel={$\omega_{1}$ / \SI{}{\radian \per \second}},
		xmin=0, xmax=30,
		ymin=-4, ymax=4,
		xtick={0,5,10,15,20,25,30},
		ytick={-4, -2,0,2,4},
		legend pos=north east,
		ymajorgrids=true,
		xmajorgrids=true,
		grid style=dashed,
		at=(one.below south),
		anchor=above north
	]
	
	\addplot[color=red,	mark=none, line width=0.4mm] table[x=time, y=data (omega1), col sep=comma] {./data/data_meas.csv};

	\addplot[color=green!60!black, mark=none, line width=0.4mm] table[x=time, y=data (omega1), col sep=comma] {./data/data_model_stiction.csv};
	
	\addplot[color=orange, dashed, mark=none, line width=0.4mm] table[x=time, y=data (omega1), col sep=comma] {./data/data_model.csv};

	\end{axis}

	\begin{axis}[
		name = three,
		width=8.4cm,
		height=3.4cm,
		xlabel={Time / \SI{}{\second}},
		ylabel={$\omega_{2}$ / \SI{}{\radian \per \second}},
		xmin=0, xmax=30,
		ymin=-3, ymax=3,
		xtick={0,5,10,15,20,25,30},
		ytick={-3, -2,-1,0,1,2.0,3},
		xmin=0, xmax=30,
		ymin=-4, ymax=4,
		xtick={0,5,10,15,20,25,30},
		ytick={-4, -2,0,2,4},	
		legend pos=north east,
		ymajorgrids=true,
		xmajorgrids=true,
		grid style=dashed,
		at=(two.below south),
		anchor=above north
		]
		\addplot[color=red,	mark=none, line width=0.4mm] table[x=time, y=data (omega2), col sep=comma] {./data/data_meas.csv}; \label{validation:p1}
	
		\addplot[color=orange,mark=none, line width=0.4mm] table[x=time, y=data (omega2), col sep=comma] {./data/data_model.csv};  \label{validation:p2}
		
		\addplot[color=green!60!black, mark=none, line width=0.4mm] table[x=time, y=data (omega2), col sep=comma] {./data/data_model_stiction.csv};  \label{validation:p3}
	
	\end{axis}

\end{tikzpicture}
	\caption{Validation of the inertia wheel pendulum: Comparison of the  mathematical models considering static friction (\ref{validation:p3}) and non static friction (\ref{validation:p2}) with measurements from the laboratory demonstrator (\ref{validation:p1}).}
	\label{fig:validation}
	\end{center}
\end{figure}

\subsection{Observer performance}

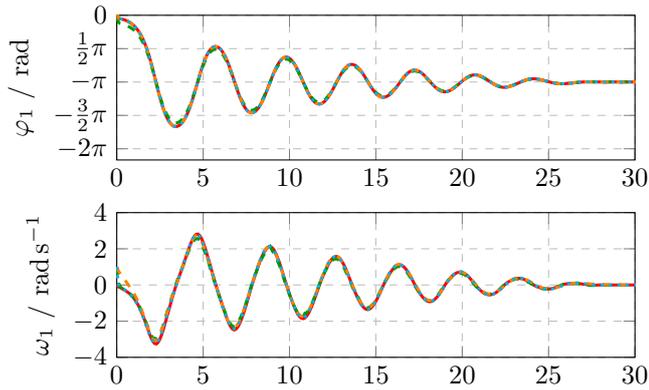
\begin{figure}[t]
	\begin{center}

\tikzsetnextfilename{\currfilebase}
\begin{tikzpicture}
	
	\begin{axis}[
	name = ov1,
	width=8.4cm,
	height=3.5cm,
	ylabel={$\varphi_{1}$ / \SI{}{\radian}},
	xmin=0, xmax=30,
	ymin=-6.81, ymax=0,
	xtick={0,5,10,15,20,25,30},
	ytick = {-2*pi,-1.5*pi,-pi,-0.5*pi,0},
	yticklabels = {$- 2\pi$, $-\frac{3}{2}\pi$,$- \pi$,$\frac{ 1}{2}\pi$, $ 0$},
	legend pos=north west,
	ymajorgrids=true,
	xmajorgrids=true,
	grid style=dashed,
	anchor=above north
	]

	\addplot+[color=red,	mark=none, line width=0.4mm, each nth point=\nthpoint] table[x=time, y=data (phi1), col sep=comma] {./data/data_meas.csv}; 
	
	\addplot+[color=green!60!black,dashed, mark=none, line width=0.4mm, each nth point=\nthpoint] table[x=time, y expr ={ \thisrow{data (phi1)}}, col sep=semicolon] {./data/data_NLOWPED_observer.csv};
	
	\addplot+[color=orange,dashed, mark=none, line width=0.4mm, each nth point=\nthpoint] table[x=time, y expr ={ \thisrow{data (phi1)}}, col sep=semicolon] {./data/data_EKF_observer.csv};
	
	\addplot+[color=cyan,dotted, mark=none, line width=0.4mm, each nth point=\nthpoint] table[x=time, y expr ={ \thisrow{data (phi1)}}, col sep=semicolon] {./data/data_LIN_observer.csv};

	\end{axis}

	\begin{axis}[
	name = ov2,
	width=8.4cm,
	height=3.5cm,
	ylabel={$\omega_{1}$ / \SI{}{\radian \per \second}},
	xmin=0, xmax=30,
	ymin=-4, ymax=4,
	xtick={0,5,10,15,20,25,30},
	ytick={-4, -2,0,2,4},
	legend pos=north east,
	ymajorgrids=true,
	xmajorgrids=true,
	grid style=dashed,
	at=(ov1.below south),
	anchor=above north
	]
	
	\addplot+[color=red,	mark=none, line width=0.4mm, each nth point=\nthpoint] table[x=time, y=data (omega1), col sep=comma] {./data/data_meas.csv}; \label{overview:meas}
	
	\addplot+[color=green!60!black,dashed,mark=none, line width=0.4mm, each nth point=\nthpoint] table[x=time, y=data (omega1), col sep=semicolon] {./data/data_NLOWPED_observer.csv}; \label{overview:NOP}
	
	\addplot+[color=orange,dashed,mark=none, line width=0.4mm, each nth point=\nthpoint] table[x=time, y=data (omega1), col sep=semicolon] {./data/data_EKF_observer.csv};  \label{overview:EKF}
	 
	\addplot+[color=cyan,dotted, mark=none, line width=0.4mm, each nth point=\nthpoint] table[x=time, y expr ={ \thisrow{data (omega1)}}, col sep=semicolon] {./data/data_LIN_observer.csv}; \label{overview:NOL}
	
	\end{axis}

\end{tikzpicture}
		\caption{Overview of the three observers: Comparison of the EKF (\ref{overview:EKF}), $\text{NO}_\text{L}$ (\ref{overview:NOL}), and $\text{NO}_\text{P}$ (\ref{overview:NOP}) estimates $\varphi_{1}$ and $\omega_{1}$ with measurements from the laboratory demonstrator (\ref{overview:meas}).}
		\label{fig:observer_overview}
	\end{center}
\end{figure}

\begin{figure}[h]
	\begin{center}

\tikzsetnextfilename{\currfilebase}
\begin{tikzpicture}
	\begin{axis}[
		name = p1,
		thick,
		no markers,
		width=8.4cm,
		height=3.5cm,
		ylabel={$\omega_{2}$ / \SI{}{\radian \per \second}},
		xmin=0, xmax=30,
		ymin=-3, ymax=3,
		xtick={0,5,10,15,20,25,30},
		ytick={-3, -2,-1,0,1,2.0,3},
		xmin=0, xmax=30,
		ymin=-4, ymax=4,
		xtick={0,5,10,15,20,25,30},
		ytick={-4, -2,0,2,4},	
		legend pos=north east,
		ymajorgrids=true,
		xmajorgrids=true,
		grid style=dashed,
		axis on top=true,
		anchor=above north
	]
	
	\addplot+[color=red,mark=none, line width=0.4mm] table[x=time, y=data (omega2), col sep=comma, each nth point=\nthpoint] {./data/data_meas.csv}; 
	
	\addlegendentry{EKF};
	
	\addplot+[color=orange,mark=none, line width=0.4mm] table[x=time, y=data (omega2), col sep=semicolon, each nth point=\nthpoint] {./data/data_EKF_observer.csv};
	
	\addplot+[name path = A,color=red,line width = 0.00001mm] table[x = time, y expr ={ \thisrow{state}==1?4:-4}, col sep=semicolon]{./data/data_EKF_observer.csv};
	
	\addplot+[name path = B,color=green,line width = 0.00001mm] table[x = time, y expr ={ \thisrow{state}==1?-4:4}, col sep=semicolon]{./data/data_EKF_observer.csv};
	
	\addplot[name path = G,mark=none, line width = 0.00001mm] coordinates {(0,-4) (30,-4)};
	
	\addplot+[red!10] fill between[of=A and G];
	\addplot+[green!10] fill between[of=B and G];
	\end{axis}
	
	\begin{axis}[
		name = p2,
		thick,
		no markers,
		width=8.4cm,
		height=3.5cm,
		ylabel={$\omega_{2}$ / \SI{}{\radian \per \second}},
		xmin=0, xmax=30,
		ymin=-3, ymax=3,
		xtick={0,5,10,15,20,25,30},
		ytick={-3, -2,-1,0,1,2.0,3},
		xmin=0, xmax=30,
		ymin=-4, ymax=4,
		xtick={0,5,10,15,20,25,30},
		ytick={-4, -2,0,2,4},	
		legend pos=north east,
		ymajorgrids=true,
		xmajorgrids=true,
		grid style=dashed,
		axis on top=true,
		at=(p1.below south),
		anchor=above north
	]
	
	\addplot+[color=red,	mark=none, line width=0.4mm, each nth point=\nthpoint] table[x=time, y=data (omega2), col sep=comma] {./data/data_meas.csv};
	
	\addplot+[color=orange,mark=none, line width=0.4mm, each nth point=\nthpoint] table[x=time, y=data (omega2), col sep=semicolon] {./data/data_LIN_observer.csv};
		
	\addplot+[name path = A,color=red,line width = 0.00001mm] table[x = time, y expr ={ \thisrow{state}==1?4:-4}, col sep=semicolon]{./data/data_LIN_observer.csv};
	
	\addlegendentry{$\text{NO}_\text{L}$}
	
	\addplot+[name path = B,color=green,line width = 0.00001mm] table[x = time, y expr ={ \thisrow{state}==1?-4:4}, col sep=semicolon]{./data/data_LIN_observer.csv};
	
	\addplot[name path = G,mark=none, line width = 0.00001mm] coordinates {(0,-4) (30,-4)};
	
	\addplot+[red!10] fill between[of=A and G];
	\addplot+[green!10] fill between[of=B and G];
	
	\end{axis}
	
	\begin{axis}[
		name = p3,
		thick,
		no markers,
		width=8.4cm,
		height=3.5cm,
		ylabel={$\omega_{2}$ / \SI{}{\radian \per \second}},
		xmin=0, xmax=30,
		ymin=-3, ymax=3,
		xtick={0,5,10,15,20,25,30},
		ytick={-3, -2,-1,0,1,2.0,3},
		xmin=0, xmax=30,
		ymin=-4, ymax=4,
		xtick={0,5,10,15,20,25,30},
		ytick={-4, -2,0,2,4},	
		legend pos=north east,
		ymajorgrids=true,
		xmajorgrids=true,
		grid style=dashed,
		axis on top=true,
		at=(p2.below south),
		anchor=above north
	]
	
	\addplot[color=red,	mark=none, line width=0.4mm, each nth point=\nthpoint] table[x=time, y=data (omega2), col sep=comma] {./data/data_meas.csv};   \label{dropdown:meas}

	\addlegendentry{$\text{NO}_\text{P}$}
	
	\addplot[orange,mark=none, line width=0.4mm, each nth point=\nthpoint] table[x=time, y=data (omega2), col sep=semicolon] {./data/data_NLOWPED_observer.csv}; \label{dropdown:observer} 
	
	\addplot+[name path = A,color=red,line width = 0.00001mm] table[x = time, y expr ={ \thisrow{state}==1?4:-4}, col sep=semicolon]{./data/data_NLOWPED_observer.csv}; 
	\addplot+[name path = B,color=green,line width = 0.00001mm] table[x = time, y expr ={ \thisrow{state}==1?-4:4}, col sep=semicolon]{./data/data_NLOWPED_observer.csv};

	\addplot[name path = G,mark=none, line width = 0.00001mm] coordinates {(0,-4) (30,-4)};
	
	\addplot+[red!10] fill between[of=A and G]; \label{red_area}
	\addplot+[green!10] fill between[of=B and G]; \label{green_area}

	\end{axis}
\end{tikzpicture}
		\caption{Estimation of the wheel angular velocity: Comparison of the proposed observer with probabilistic model selection (\ref{dropdown:observer}) and measurements from the laboratory demonstrator (\ref{dropdown:meas}).}
		\label{fig:observer_omega2}
	\end{center}
\end{figure}
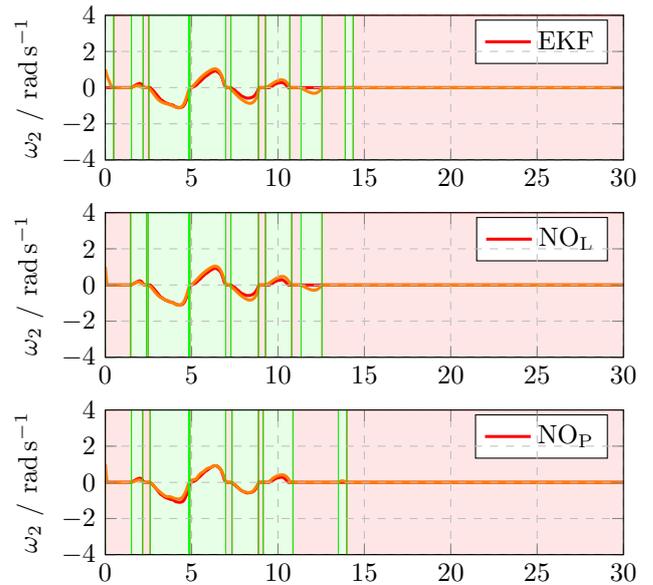

The three proposed observers with probabilistic model switching are presented in the context of the drop-down experiment, i.\,e., free swinging down of the inertia pendulum to the lower stable rest position.
The observer parameters were not optimized, since a direct comparison was not intended, rather the principle suitability should be demonstrated.
Apart from the introductory overview, see Fig. \ref{fig:observer_overview}, the observer behavior with respect to the angular velocity $\omega_{2}$ is investigated in detail, since the pendulum angle is used as a measurand and the static friction particularly affects the wheel. 
The transient behavior of the estimates results from the intentionally miss initialized observers. 
The experiment starts in the state $x_0 = [0,0,0]^T$, while the observers are initialized with $\hat{x}_0 = [\SI{-\pi/10 }{\radian}, \SI{1}{\radian \per \second},\SI{1}{\radian \per \second}]^T$. The covariances of the Gaussian distributed quantities for the EKF and the $\text{NO}_\text{L}$, see Section \ref{sec:observer}, are set to $P = \text{diag}(0.00165, 0.01, 0.1)$, $Q = \text{diag}(0, 0.01, 0.1)$ and $R = 0.001$.
The nonlinear observer with linear error dynamics has been implemented as sampled-data system (\ref{equ:nlowled_sampled_data}).
Designed as a time-invariant linear-quadratic estimator, the correction term of the error dynamics is obtained as Kalman gain.
Assuming that the measurement is not disturbed, i.\,e., $R = 0$, the solution of the algebraic Riccati equation yields observer parameter that lead to a dead-beat behavior. The parameter of nonlinear observer with passive error dynamics are set to $\alpha = 10$ and $\beta = 5$. 
{The observer estimates in Fig. \ref{fig:observer_overview} are nearly identical, except for the transient responses.}
The estimates of the wheel angular velocity $\hat{\omega}_2$ are presented in Fig. \ref{fig:observer_omega2}. 
The results of the observers are divided into three separate graphs to illustrate the change in adhesion of the wheel. 
In every estimation step, the observers first determine the probability that the system is in stiction and not in stiction given the latest measurement.
The most likely model is subsequently applied for the estimation of the state. 
In Fig. \ref{fig:observer_omega2}, the current model used for estimation is indicated by the background color, with a red area (\ref{red_area}) marking the sticking model and a green area (\ref{green_area}) marking the no-sticking model.
{All three observers provide satisfactory results, especially with respect to their combination with the probabilistic model selection.}

\section{Conclusion}

In this contribution, three observers for the nonlinear model of inertia wheel pendulum with switching stiction of the wheel have been presented. Beside an extended Kalman Filter, a nonlinear observer with linear as well as a nonlinear observer with passive error dynamics are proposed. As the considered system is described by two differential equations, representing a sticking and a non-sticking model, the switching in the observer is implemented by means of a probabilistic model selection based on Bayes factor.
The performance of the observers is verified on a laboratory demonstrator.

\bibliography{ifacconf}             
                                                   







\end{document}